\documentclass[11pt,reqno,a4paper]{amsart}
\usepackage{amsmath}
\usepackage{amsfonts}
\usepackage{amssymb}
\usepackage{breqn}
\usepackage[latin1]{inputenc}
\usepackage{graphicx}

 \newtheorem{thm}{Theorem}[section]
 \newtheorem{cor}[thm]{Corollary}
 
 \newtheorem{prop}[thm]{Proposition}
 \theoremstyle{definition}
 \newtheorem{defn}[thm]{Definition}
 
 \theoremstyle{remark}
 \newtheorem{rem}[thm]{Remark}
 \newtheorem{ex}{Example}
 \numberwithin{equation}{section}

\def\Z{{\mathbb Z}}
\def\N{{\mathbb N}}
\def\R{{\mathbb R}}
\def\C{{\mathbb C}}

\def\codiff{\mbox{codiff}}
\def\cosum{\mbox{cosum}}

\def\codiff{\mbox{codiff}}
\title[Integral representations for special functions]
 {Positive-definiteness and integral representations for special functions}
\author[J. Buescu]{J. Buescu}
\address{%
Dep. Matem\'{a}tica\\ FCUL and CMAFCIO \\ Portugal}

\email{jsbuescu@fc.ul.pt}

\thanks{The first author acknowledges partial support by  Fundação para a Ciência e
Tecnologia, UID/MAT/04561/2013. }
\author{A. C. Paixão}
\address{Área Departamental de Matemática\\ ISEL\\ Portugal }
\email{apaixao@dem.isel.ipl.pt}



\date{2017-01-29}

\begin{document}

\begin{abstract}
We characterize a holomorphic positive definite function $f$ defined on a horizontal strip of the complex 
plane as the Fourier-Laplace transform of a unique exponentially finite measure on $\R$. 
The classical theorems of Bochner on positive definite functions and of Widder on exponentially convex functions 
become special cases of this characterization: they are respectively the real and  pure imaginary sections of
the  complex integral representation. 
We apply this integral representation to special cases, including the $\Gamma$, $\zeta$ and Bessel functions, and 
construct explicitly the corresponding measures, thus providing new insight into the nature of complex 
positive and co-positive definite functions: in the case of the zeta function this process leads to a 
new proof  of an integral representation on the critical strip. 
\end{abstract}

\maketitle

%
%
%
%
%
%

\section{Introduction}


This is an overview 
of some fundamental representations of positive definite functions and their connections to regularity 
and extension properties of these functions.
In this paper we derive consequences of these connections for the characterization of 
complex variable positive definite functions on horizontal strips of the complex plane. 
This point of view is used in \S \ref{sec_moments_MGF_CF} as a means to establish connections 
between the concepts of characteristic function, moment generating function and complex positive definiteness.

The paper is structured as follows. The first section will be dedicated to the presentation of basic definitions and relevant properties 
of the different notions of positive definiteness involved in the subsequent discussion. The conditions under which a complex variable
positive definite function defined on a general (horizontal or vertical) 
strip admits an integral transform representation are next explored. 
This point of view is then applied to specific cases, including the special functions $\Gamma (z)$, $\zeta(z)$
and Bessel functions of the first kind. 
As it turns out, $\Gamma$ is  alternately co-positive and co-negative definite on neighbouring vertical strips of the complex plane,
 $\zeta$ is co-positive definite on the semiplane $|z| > 1$, and the normalized Bessel functions are co-positive definite 
on the entire complex plane. For these functions we compute explicitly the 
exponentially finite measures relative to which the integral representation holds, and as a corollary 
we derive a new proof of a known integral representation of the $\zeta$ function on the critical strip.

This paper is partly expository in nature. The main theorems have been published elsewhere; 
proofs are omitted and, unless otherwise stated, may be found in references 
\cite{bp_DIE, bp_JIEA, bp_LAA, bp_pdmdrki, bp_diff_PDF, bp_RC_PDF, bp_CVPDF, bp_CPDFS, bp_PROP} by the authors.
The results involving special functions, namely the $\Gamma$, $\zeta$ and Bessel functions, are original.


\section{Real variable positive definite functions}
\label{subsec_RVPDF}

A brief review of the properties of real variable positive definite functions 
relevant for our purposes will be performed below.

A function $f: \R \to \C$ is positive definite if 
\begin{equation}
\label{eq_posdef}
 \sum_{k,l = 1}^n f(x_k - x_l) \xi_k \, \overline{\xi_l} \geq 0 
\end{equation}
for every choice of $x_1, \ldots x_n \in \R$ and $\xi_1, \ldots, \xi_n \in \C$; that is, if every 
matrix $ [f(x_k - x_l)]_{k,l = 1}^n$ is positive semidefinite. 

Basic properties of positive definite functions include the following well-known facts: 
\begin{enumerate}
\item[1.] $f(0) \geq 0$ and $f(-x) = \overline{f(x)} \ \  \forall x \in \R$,
\item[2.] $|f(x)| \, \leq \, f(0) \ \  \forall x \in \R.$
\end{enumerate}
These follow immediately form \eqref{eq_posdef} by consideration of order 1 and 2 matrices; references 
to these properties may be found in  \cite{mat}, \cite{ste}. 
We stress that these are purely algebraic properties that hold under the single hypothesis of positive definiteness
and for which no form of regularity is necessary. Under the further hypothesis of continuity, positive definite functions 
satisfy the following classical representation theorem.

\begin{thm}[Bochner]
\label{thm_Bochner}
A continuous function $f : \R \to \C$ is positive definite if and only if it is 
the Fourier transform of a unique, finite, non-negative measure $\mu$ on $\R$, that is
\begin{equation}
\label{eq_Bochner} 
f(x) = \int_{-\infty}^{+\infty} e^{itx} \, d\mu(t). 
\end{equation}
\end{thm}

The Bochner representation \eqref{eq_Bochner} is often the most expedite way to prove that a function is positive definite. 
Indeed, the most obvious class of positive definite functions are characteristic functions, which 
arise naturally in probability theory as Fourier transforms of probability measures.  Specific examples follow immediately: 
functions like $\cos (ax), \, e^{-|x|}, \, \frac{1}{1+x^2}, \, e^{-x^2}$ are Fourier transforms of
probability measures; consequently, 
from the Bochner representation theorem, they are positive definite.

A first case of propagation of regularity concerns the consequences of certain regularity assumptions for a real variable positive 
definite function on a neighborhood of the origin.

\begin{thm}
\label{thm_propagation1}
Let  $\phi$ be a positive definite function. Suppose $\phi$ is of class $C^{2k}$ on a neighborhood of the origin. 
Then $\phi$ is $C^{2k}(\R)$.
\end{thm}

This is a folklore theorem for which several proofs exist. The classical proof relies in standard properties of the
Fourier transform and may be found, for instance, in \cite{don}. Constructive approaches based on the algebraic properties 
arising directly from the definition may be may be found in \cite{bp_diff_PDF} or \cite{bp_PROP}. 
%
%
%
%
%
%

\section{Holomorphic extensions}
\label{subsec_HolExt}

We next turn to an extension result which states that it is possible to extend a 
real-analytic positive definite function as a holomorphic function to a maximal horizontal strip of the complex plane.
Priority seems difficult to establish; for proofs see \cite{luk2}, \cite{sas94, sas00, sas13}. A proof along different lines is also supplied
 in \cite{bp_CVPDF}.

\begin{rem}
\label{rem_analiticity}
Note that a function $f : \R \to \C$ is analytic in the open set $A \subset \R$ if and only if there exists a complex function $\Theta$ 
of one complex variable which is holomorphic on an open set $S \subset \C$ such that $f(x) = \Theta (x)$ for all $ x \in A$.
\end{rem}

\begin{defn}
\label{defn_strips}
A horizontal strip on the complex plane is an open set of the form
\[ S_{a,b} = \{ z \in \C : a < \Im(z) <b \}    \mbox{ for } a, \, b \in \overline{\R}.\]
Analogously, a vertical strip on the complex plane is a set of the form
\[T_{a,b} = \{ z \in \C : a < \Re(z) <b \} \mbox{ for } a, \, b \in \overline{\R}.\]
\end{defn}

\begin{thm}
\label{thm_sasvari}
Suppose 
$f: \R \to \C$ is a positive definite function which is analytic in a neighborhood of the origin. Then $f$ is analytic in $\R$
and there exist $\alpha, \beta \in (0, +\infty]$ such that $f$ extends to a function which is holomorphic in the horizontal strip 
$S_{-\alpha, \beta}$, with $\alpha$ and $\beta$ maximal with respect to this property. If $\alpha < \infty$ (resp. $\beta < \infty$),
then $-i\alpha$ (resp. $i\beta$) is a singularity of this extension, which may be defined on this strip by 
\[ \tilde{f}(z) = \int_{-\infty}^\infty \, e^{izt} \, d\mu(t), \ \ \ -\alpha < \Im(z) < \beta, \]
where $\mu$ is a finite non-negative measure on $\R$.
\end{thm}

Further reference to related results and generalizations may be found in \cite{bp_CPDFS}. 

With a view to addressing another line of holomorphic extension results with special significance to the contents of the next section, 
we now turn our attention to real variable co-positive definite functions or, as they are traditionally referred to, exponentially
convex functions.

\begin{defn}
\label{def_co-positive}
A function $G: I \to \C$, with $I \subset \R$, 
is said to be co-positive definite if, for every $n \in \N$ and every finite collection
$\{ x_j \}_{j=1}^n$ such that $x_k + x_l \in I$ for all $k, l, = 1, \ldots, n$ 
we have
\begin{equation}
\label{eq_expconv}
 \sum_{k,l = 1}^n G(x_k + x_l) \xi_k \, \overline{\xi_l} \geq 0 
\end{equation}
for every finite collection $\{ \xi_j \}_{j=1}^n \subset \C$.
\end{defn}

Co-positive definite functions are most usually defined in the case where $I$ is an interval, in which case they are known in the
literature as {\em exponentially convex functions\/}. In the real-variable context, these functions
have been introduced by Bernstein in 1929 \cite{ber} and have since been studied by 
Widder \cite{wid}, Devinatz \cite{dev2}, \cite{dev3} and other authors.

A number of properties may be established for co-positive definite functions as a direct algebraic consequence  of
definition \ref{def_co-positive}. Much more, however, can be said if we assume continuity of $G$ in an open
interval.
Widder's theorem \cite{wid}, which we rephrase in a more convenient form, asserts that the co-positive definite function 
$G$ admits the following integral representation.

\begin{thm}(Widder)
\label{thm_widder}
Let $I = (a,b)$. A function $G : I \to \R$ may be represented in the form 
\begin{equation}
\label{eq_Widder}
G(y) = \int_{-\infty}^{+\infty} e^{-yt} \, d\mu_I(t), \ \ \ y \in (a,b) 
\end{equation}
where $\mu_I$ is a non-negative Lebesgue-Stieltjes measure, if and only if $G$ is continuous and co-positive definite on $I$. 
\end{thm}

From Theorem \ref{thm_widder} we may deduce the following 

\begin{cor}
\label{cor_complexPD}
 Suppose $G(x)$ is a continuous co-positive definite function on the interval $I=(a,b)$. Then $G$ has a 
 holomorphic extension $g(z)$ to the vertical strip $T_{a,b}$
 defined by 
 \begin{equation}
\label{eq_Widder2}
  g(z)=\int_{-\infty}^{+\infty}e^{-zt}d\mu_I(t)
 \end{equation}
 that is co-positive definite in $T_{a,b}$. Furthermore, we have
\begin{equation}
\label{eq_derivatives_codp}
\frac{d^n}{dz^n} G(z) = \int_{-\infty}^{+\infty} (-t)^n \, e^{-zt} d\mu_I.
\end{equation}
\end{cor}

The contents and interplay of theorems  \ref{thm_sasvari} and \ref{thm_widder} may be viewed as the natural
motivation for the definition of positive and co-positive definite functions of a complex variable. 
This is the main concern of the next section, along with the study of corresponding consequences for the integral
representation of these functions.

%
%
%
%
%
%

\section{Complex variable positive and co-positive definite functions}
\label{sec_complex_PDF_coPDF}

 We begin with a short introduction to complex variable positive definite functions. A more comprehensive discussion may be found in 
\cite{bp_RC_PDF}, \cite{bp_CVPDF}.

\begin{defn}
\label{def_codifference}
A set $S \subset \C$ (resp. $T \subset \C$) is called a {\em codifference set\/} (resp. {\em cosum set\/})
if there exists a set $\Omega \subset \C$ such that $S$ may be written as
\[ S= \Omega - \overline{\Omega} \equiv  \{ z \in \C : \exists \  z_1, z_2 \in \Omega:
z = z_1- \overline{z_2} \}, 
\]
resp.
\[ T= \Omega + \overline{\Omega} \equiv  \{ z \in \C : \exists \  z_1, z_2 \in \Omega:
z = z_1 + \overline{z_2} \}.
\]
We write $S= \codiff(\Omega)$ (resp. $T = \cosum(\Omega)$). 
\end{defn}

Basic but very useful examples of codifference and cosum sets are, respectively, horizontal and vertical
strips $S_{a,b}$ and $T_{a,b}$ as in definition \ref{defn_strips}.

\begin{defn}
\label{def_PD_coPD_in_C} Let $\Lambda \subset \C$. 
A function $f : \Lambda \to  \C$ is said to be {\em positive definite in $S \subset \Lambda$\/} 
(resp. {\em co-positive definite in $T \subset \Lambda$\/}) if, for every $n \in \N$ 
and every finite collection $\{z_j \}_{j=1}^n$ such that $z_k- \overline{z_l} \in S$ for all $k, \, 
l = 1, \ldots, n$ 
(resp. $z_k + \overline{z_l} \in S$ for all $k, \, l = 1, \ldots, n$), we have
\[\sum_{k, l=1}^n   f(z_k-\overline{z_l})\, \xi_k \,
\overline{\xi_l} \geq 0,  \]
resp. 
\[\sum_{k, l=1}^n   f(z_k+\overline{z_l})\, \xi_k \,
\overline{\xi_l} \geq 0,  \]
for every collection $\{ \xi_j \}_{j=1}^n \subset \C$, that is,  if every 
matrix $ [f(z_k - \overline{z_l})]_{k,l = 1}^n$ (resp. $ [f(z_k + \overline{z_l})]_{k,l = 1}^n$ )
is positive semidefinite.  
\end{defn}

\begin{rem}
\label{rem_horiz_imag_axis}
We observe that, for a positive definite function $f$ on a horizontal strip $S_{a,b}$, the real variable function $F_y(x) = f(x+iy)$ for fixed
$y \in (a,b)$ verifies the condition $\sum_{k,l= 1}^n F_y(x_k - x_l) \, \xi_k \overline{\xi_l} \geq 0$ for every collection $\{x_j\}_{j=1}^n \subset \R$,
every collection $\{\xi_j\}_{j=1}^n \subset \C$ and every positive integer $n$, and is thus a  positive definite
real-variable function. On the other hand, the function $G(y) = f(iy)$, for $y \in (a,b)$, satisfies the inequalities 
$\sum_{k,l= 1}^n G(x_k + x_l) \, \xi_k \overline{\xi_l} \geq 0$. 
We recognize this property as the real-variable co-positive definiteness condition. 
\end{rem}

These definitions of cosum set and co-positive definite function are readily seen to be closely related to those of codifference set
and positive definite function, respectively. More precisely, the following is true.

\begin{prop}
\label{prop_codiff_cosum} We have:
\begin{enumerate}
\item[(i)] $T \subset \C$ is a cosum set if and only if $S = iT$ is a codifference set.
\item[(ii)] $f(z)$ is positive definite in a set $S \subset \C$ if and only if $g(z) = f(iz)$ 
 is co-positive definite in $T = -iS$.
\end{enumerate}
\end{prop}

It is immediate to show that $S_{a,b}$ is a horizontal strip if and only $-i \, T_{a,b}$ is a vertical strip.  
Note also that Definition \ref{def_PD_coPD_in_C} does not require any regularity on the function $f$. Even in this context, 
a number of properties of positive definite functions may be established as purely algebraic consequences of definition 
\ref{def_PD_coPD_in_C}, as the next proposition shows for positive definite functions.

\begin{prop}
\label{prop_basic_properties}
Let $f$ be a positive definite function on a set $S \subset \C$. Then the following hold.
\begin{enumerate}
\item[(1)] If $x, y \in \R$ and $\pm x + i y \in S$, then $f(-x+iy) = \overline{f(x+iy)}$. In particular, taking $x= 0$ implies that $f(iy)$ is real
and non-negative.
\item[(2)] If $x, y, \beta \in \R$ and $\pm x + i y, \, i(y \pm \beta) \in S$, then 
\[ |f(x+iy)|^2 \, \leq \, f(i(y-\beta)) f(i(y+\beta)). \]
\item[(3)] If $f$ is positive definite, then so is $\overline{f}$.
\item[(4)] If $f_1, \ldots, f_n$ are positive definite functions and $c_i \geq 0$, then so is $\sum_{k=1}^n c_k f_k$. 
\item[(5)] If $\{f_n \}_{n \in \N}$ is a pointwise convergent sequence of positive definite functions, 
then the limit $f(z) = {\displaystyle \lim_{n \to \infty} f_n(z)}$ is a positive definite function.
\item[(6)] The product of positive definite functions is positive definite.
\end{enumerate}
\end{prop}

Corresponding properties for co-positive definite functions immediately follow from proposition \ref{prop_codiff_cosum}.

We close this section with a reference to non-trivial examples of functions defined on strips of the complex plane satisfying some kind of 
positive or negative definiteness condition. The following definition will be useful for this purpose.

\begin{defn}
\label{def_DN_co_DN}
Let $\Lambda \subset \C$.  A function $f : \Lambda \to  \C$ is said to be
\begin{enumerate}
\item[i)] negative definite  in $\Lambda$ if $g = -f$ is positive definite in $\Lambda$;
\item[ii)] co-negative definite  in $\Lambda$ if $g = -f$ is co-positive definite  in $\Lambda$. 
\end{enumerate}
\end{defn}

The following examples illustrate the richness of behaviour of complex-variable positive definite functions.

\begin{ex} The function $f(z) = \frac{1}{1+z^2}$ is positive definite in the horizontal strip $S_{-1,1}$ and negative defini\-te in
the half-planes  $S_{-\infty,-1}$ and $S_{1, \infty}$.  Let $j$ and $n$ be positive integers and define $a_j = j$
for $j \leq n$, $a_{n+1} = \infty$.
Since the product of positive definite functions is positive definite, it follows that 
\[ f_n(z) = \prod_{j=1}^n  \frac{1}{(1-i\frac{z}{a_j}) (1+i\frac{z}{a_j})}  \]
is positive definite in the central horizontal strip $S_{-1, 1}$, alternating between negative and positive definite
in neighboring horizontal strips $S_{a_j, a_{j+1}}$ and $S_{-a_{j+1}, -a_j}$ for $j=1, \ldots, n$.
\begin{figure}[h]
    \centering
    \includegraphics[width=8cm]{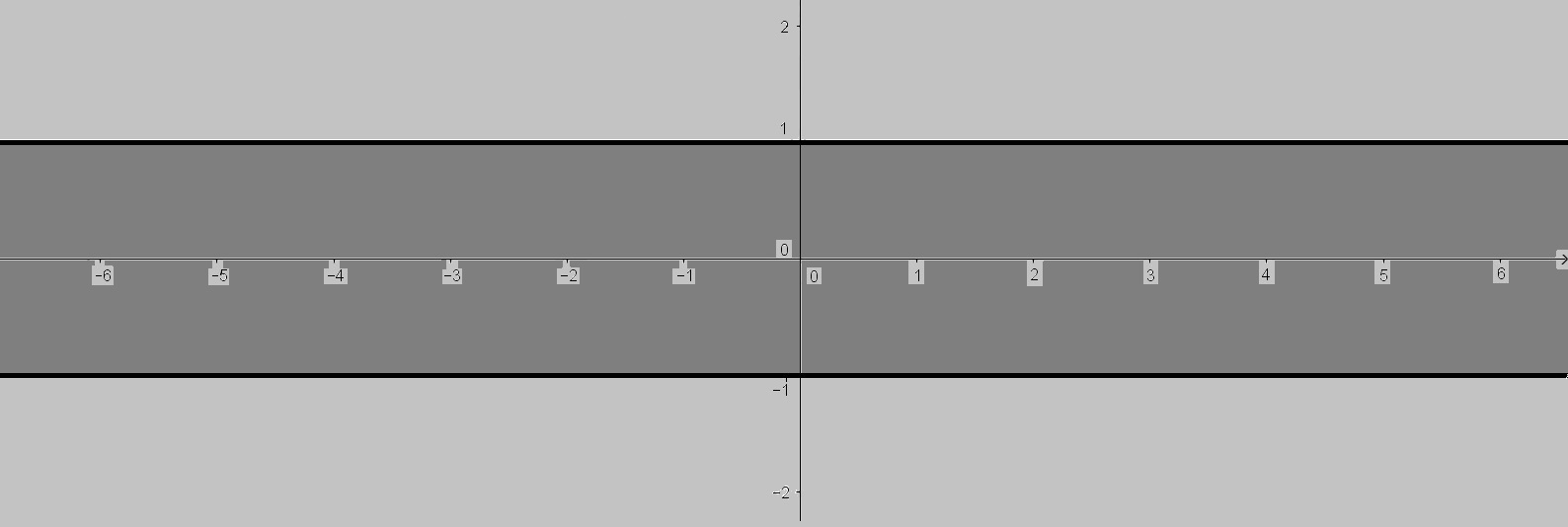}
    \caption{$f(z) = \frac{1}{1+z^2}$ is positive definite in the horizontal strip $S_{-1,1}$ and negative definite in the complementary half-planes.}
    \label{fig_exemplo_1}
\end{figure}

\end{ex}

\begin{ex} The function $f(z) = \frac{1}{\cosh (\frac{\pi z}{2})}$ is positive definite in any 
horizontal strip of the form  $S_{4n-1, 4n+1}$ and negative definite  in any 
horizontal strip   $S_{4n+1, 4n+3}$ for $n \in \Z$.
In fact, the  well-known integral representation (see e.g. \cite{kos})
\begin{equation}
\label{formula_cosh}
f(z) = \frac{1}{\cosh \left( \frac{\pi z}{2} \right)} = \frac{1}{\pi}\int_{-\infty}^\infty \frac{1}{\cosh t} e^{-izt} \, dt
\end{equation}
implies that it is positive definite on the strip $S_{-1,1}$. It is not difficult to show directly that $f$ is negative definite on 
$S_{1,3}$. Since $f$ is periodic of period $4i$, it alternates between being positive definite on horizontal strips $S_{4n-1, 4n+1}$ and
negative definite on horizontal strips $S_{4n+1, 4n+3}$.
\begin{figure}[ht]
    \centering
    \includegraphics[width=8cm]{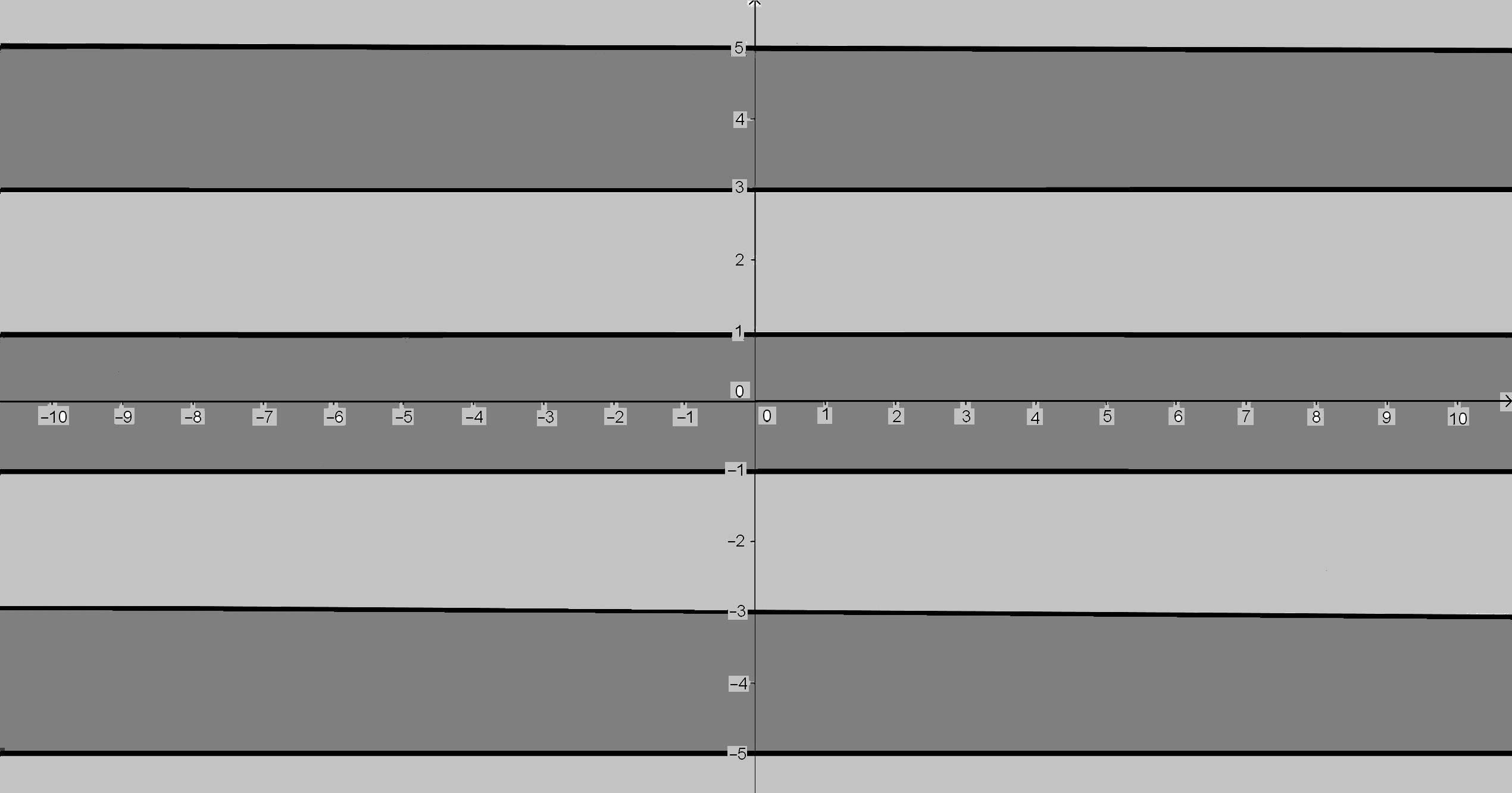}
    \caption{$f(z)= \frac{1}{\cosh (\frac{\pi z}{2})}$ is positive definite in the dark horizontal strips and negative definite on the light horizontal strips. }
    \label{fig_exemplo_2}
\end{figure}

\end{ex}

\begin{ex}
\label{ex_Gamma}
The $\Gamma$ function, defined by 
\[ \Gamma(z) = \int_0^{+\infty} x^{z-1} e^{-x} \, dx \]
for $\Re(z) > 0$, 
is easily seen to satisfy definition \ref{def_PD_coPD_in_C}, so it is co-positive definite in that semi-plane. The $\Gamma$ function 
has a simple pole at $z=0$ and at every negative integer. In each vertical strip $T_{-n-1,-n}$, for  integer $n \geq 0$, the $\Gamma$ function
satisfies the Cauchy-Saalschütz formula (see e.g. \cite{whi-wat}, pg. 244)
\[ \Gamma(z) = \int_0^{+\infty} x^{z-1} \left( e^{-x} - \sum_{m=0}^n   (-1)^m \frac{x^m}{m!}  \right)    \, dx,  \]
from which a simple calculation shows that it is co-negative definite on the vertical strip $T_{-1,0}$, alternating 
from co-negative definite to co-positive definite and back in neighbouring vertical strips. Corresponding statements are valid 
for positive-definiteness of $\Gamma(-iz)$ on horizontal strips in view of Proposition \ref{prop_codiff_cosum}.
\begin{figure}[ht]
    \centering
    \includegraphics[width=8cm]{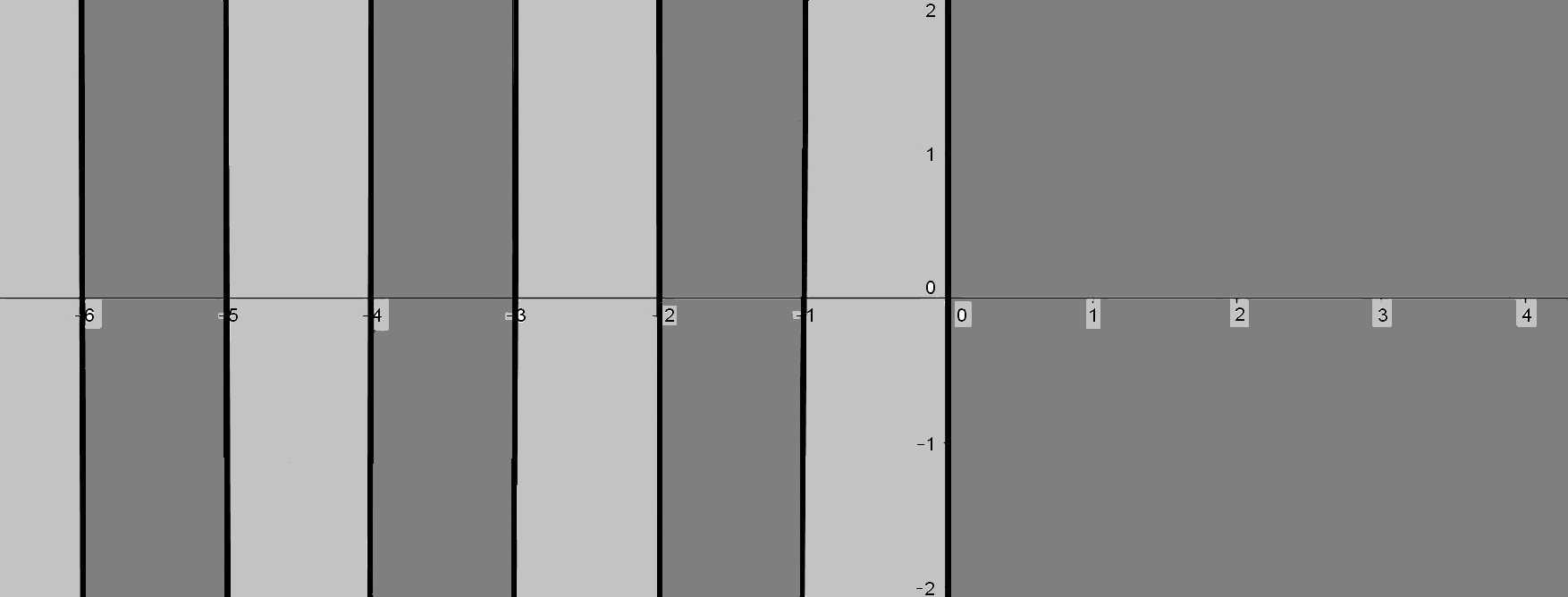}
    \caption{The $\Gamma$ function is co-positive definite on the half-plane $\Re(z) > 0$ and in the dark shaded vertical strips, and 
		co-negative definite on the light vertical strips.}
    \label{fig_exemplo_3}
\end{figure}

\end{ex}

\begin{ex}
\label{ex_Bessel}
Consider the real-variable Bessel function of the first kind $J_\alpha(x)$. 
For $\alpha > -\frac{1}{2}$ the normalized Bessel function 
$\mathcal{J}_\alpha(x)$ is defined by
\[ \mathcal{J}_\alpha (x) = 2^\alpha \Gamma(\alpha+1) x^{-\alpha} J_\alpha(x) \]
extends to a real-analytic in $\R$ with $\mathcal{J}_\alpha(0) = 1$ by, see \cite{ask, bar, neu, sel}. 
This function satisfies the integral representation 
\[ \mathcal{J}_\alpha(x) = \frac{1}{2^{\alpha-1} \sqrt{\pi} \Gamma(\alpha+1/2)} \int_0^1 (1-t^2)^{\alpha -1/2} \, \cos(xt) \, dt, \]
see e.g. Watson \cite{wat}.
This representation admits the natural extension to the complex plane 
\[  \mathcal{J}_\alpha(z) = \frac{1}{2^{\alpha-1} \sqrt{\pi} \Gamma(\alpha+1/2)} \int_0^1 (1-t^2)^{\alpha-1/2} \, \cos(zt) \, dt. \]
This integral may be rewritten in the form of a Fourier-Laplace transform
(see example \ref{ex_Bessel_1st_kind}), and application of theorem 
\ref{thm_characterization1} below shows that $\mathcal{J}_\alpha$ is an entire 
complex-variable positive definite function. 
This implies, for integer $\alpha \geq 0$, that 
\[ J_\alpha(z) = \frac{1}{2^\alpha \Gamma(\alpha+1)} z^\alpha \mathcal{J}_\alpha(z) \]
defines an entire extension $J_\alpha$ to the complex plane.
If $\alpha = 0$, in particular, we have $\mathcal{J}_0(z) = \mathcal{J}_0(z)$ and conclude that 
$J_0$ is an entire positive definite function.
\end{ex}

\begin{ex}
\label{ex_Zeta}
Riemann's zeta function $\zeta(z)$ is defined by
$ \zeta(z) = \sum_{n=1}^\infty \frac{1}{n^z}$  for $z \in T_{1,+\infty}$. Each summand in the series
 $\frac{1}{n^z} = e^{- z \ln n}$ is co-positive definite and therefore  $\zeta(z) = \sum_{n=1}^\infty \frac{1}{n^z}$ is co-positive definite
on the half-plane of convergence $T_{1,+\infty}$ by Prop. \ref{prop_basic_properties}. 
The zeta function is meromorphic and has a unique simple pole at
$z=1$. 
However, in contrast to previous examples,  $\zeta(z)$ is neither co-positive definite nor co-negative definite in any 
open cosum set in the half-plane $T_{-\infty, 1} = \{ z \in \C: \Re(z) < 1 \}$ since, otherwise, $f$ would be co-positive (resp. co-negative)
definite on some open interval $I \subset (-\infty, 1)$ (refer to \cite{bp_CVPDF} and use Prop. \ref{prop_codiff_cosum} to derive
the relevant properties of cosum sets). By Theorem \ref{thm_minimal_PD} below, $\zeta(z)$ would necessarily be co-positive
(resp. co-negative) in the whole half-plane $T_{-\infty, 1}$. The existence of the trivial zeros of the zeta function would then imply, 
by the inequalities corresponding to statement (2) in Proposition \ref{prop_basic_properties}, 
that $\zeta(z)$ would vanish identically in the whole half-plane $T_{-\infty, 1}$,
which is not the case.
\begin{figure}[ht]
    \centering
    \includegraphics[width=8cm]{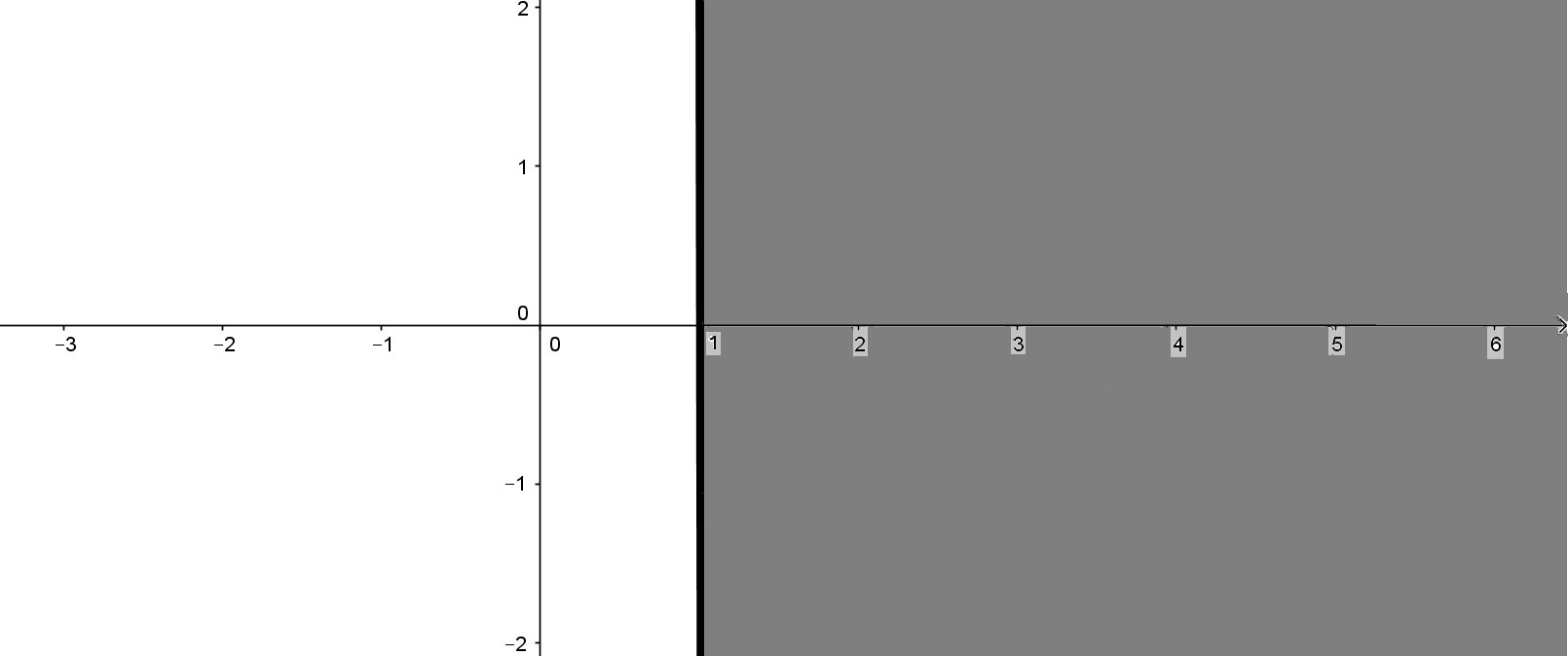}
    \caption{The $\zeta$ function is co-positive definite on the half-plane $\Re(z) > 1$. On the complementary half-plane it
		is neither co-positive nor co-negative definite.}
    \label{fig_exemplo_4}
\end{figure}
\end{ex}

\begin{ex}
\label{ex_function_ZED}
Define the function 
\begin{equation}
\label{eq_def_ZED} 
	\mathcal{Z}(z) = -\frac{\zeta(z)}{z}. 
\end{equation}
It is not difficult to show that $\mathcal{Z}$ is co-negative definite on the half-plane $\Re(z) > 1$, 
co-positive definite on the vertical strip $0 < \Re(z) < 1$ and neither co-positive neither co-negative
on the half-plane $\Re(z) <0$. 

\begin{figure}[ht]
    \centering
    \includegraphics[width=8cm]{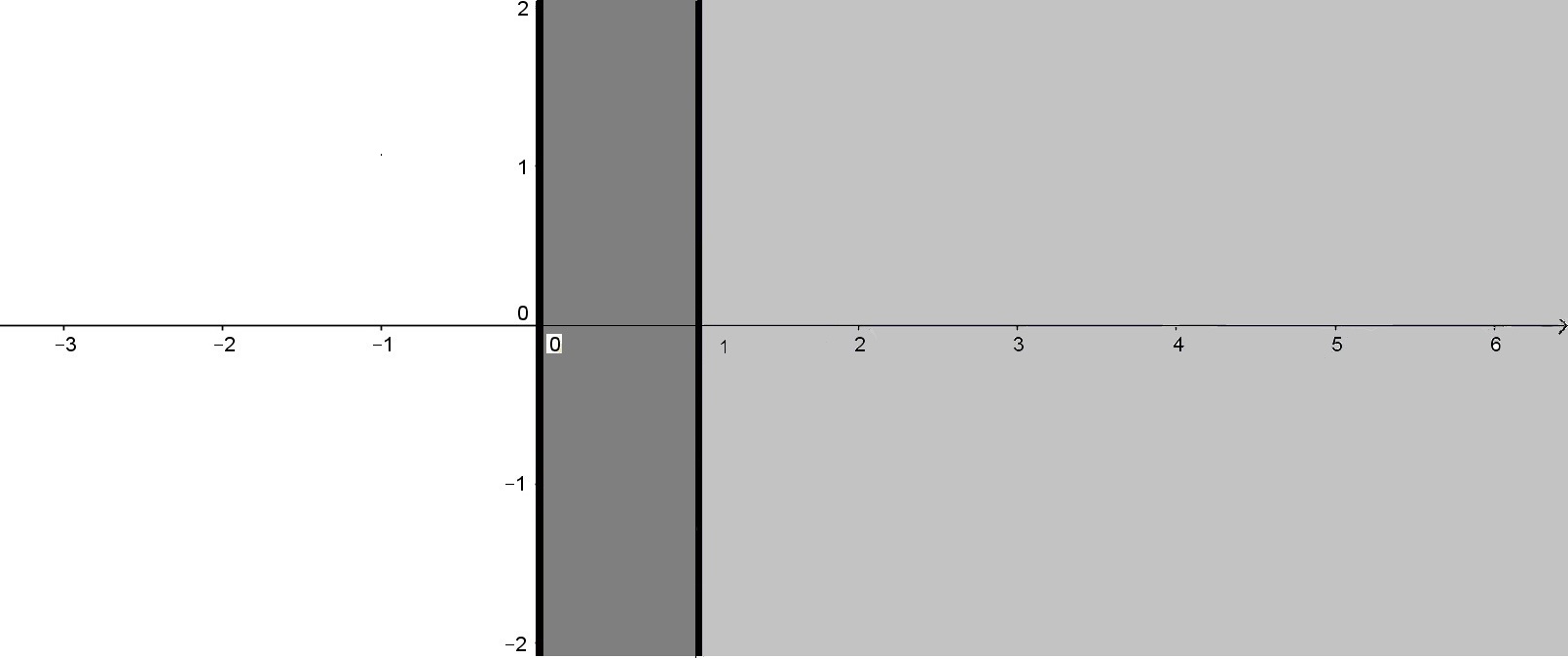}
    \caption{The $\mathcal{Z}$ function is co-negative definite on the half-plane $\Re(z) > 1$ (light grey) and 
		co-positive definite on the critical strip $0 < \Re(z) < 1$ (dark grey). On the  half-plane $\Re(z) < 0$ 
		(white) it is neither co-positive nor co-negative definite.}
    \label{fig_example_ZED}
\end{figure}

A full discussion of these properties of $\mathcal{Z}$ is deferred to 
example \ref{example_ZED}, where the interest of introducing this function will become apparent.
\end{ex}

%
%
%
%
%
%

\section{Positive definite functions as Fourier-Laplace transforms}
\label{sec_PDF_FLT}

The concept of {\em exponentially finite measure\/}, which we now introduce, will prove essential for what follows.

\begin{defn}
\label{def_exp_fin_measure}
A non-negative measure $\mu$ is said  to be {\em exponentially finite with respect to the (non-empty) interval $I$} if
$\int_{-\infty}^{+\infty} e^{-yt} \, d\mu(t)$ is finite for every $y \in I$.
\end{defn}

\begin{rem}
\label{rem_exp_fin_are_LS}
It is easily seen that an exponentially finite measure is a Lebesgue-Stieltjes measure. 
\end{rem}

Under the hypothesis of holomorphy, the notions of positive and co-positive definiteness reviewed in the previous sections 
may now be related in the following way. 

\begin{thm}
\label{thm_complex_equivalence}
Suppose $f$ is a holomorphic complex function defined on the horizontal strip 
$S_{a,b}$. Define $F_y(x) = f(x + iy)$ 
for fixed $y \in (a, b)$ and $x \in \R$, and $G(y) = f(iy)$ for 
$y \in (a,b)$. Then the following statements are equivalent: 

\begin{enumerate}
\item $f$ is a complex-variable positive definite function on $S_{a,b}$.
\item There exists a non-negative Lebesgue-Stieltjes measure $\mu$ on $\R$ such that 
\[ f(z) = \int_{-\infty}^{+\infty} e^{itz} \, d\mu(t) \ \ \ \forall z \in S_{a,b}. \]
\item $F_y(x)$ is a positive definite function of the real variable $x$ for some $y \in (a,b)$.
\item There exists a non-negative finite measure $\mu_y$ on $\R$ such that
\[ F_y(x) = \int_{-\infty}^{+\infty} e^{itx} \, d\mu_y(t) \ \ \ \forall x \in \R \]
for some $y \in (a, b)$. 
\item $G(y)$ is a co-positive definite function of the real variable $y$ on the interval $(a,b)$.
\item There exists a non-negative Lebesgue-Stieltjes measure $\mu_I$ on $\R$ such that 
\[ G(y) = \int_{-\infty}^{+\infty} e^{-yt} \, d\mu_I(t) \ \ \ \forall y \in (a,b). \]
\end{enumerate}
Furthermore, in any of these conditions we also have: 
\begin{enumerate}
\item[(i)] The measures $\mu = \mu_I$ and $\mu_y$ are uniquely determined by the function $f$.
\item[(ii)] $\mu = \mu_I$ is exponentially finite with respect to $I = (a,b)$ and $\mu_y$ is exponentially finite
with respect to $(a-y, b-y)$.
\item[(iii)] All the functions $f(z),  F_y(x)$ and $G(y)$ are infinitely differentiable and may be differentiated 
under the integral sign in their respective representations (2), (4) and (6).
\end{enumerate}
\end{thm}

Theorem \ref{thm_complex_equivalence} may be seen as a consequence of the interplay of theorems 
\ref{thm_sasvari} and \ref{thm_widder} (along with Corollary \ref{cor_complexPD})
which reveals, through consideration of the complex variable setting, the precise joint significance 
of these individual results.

We are now ready to present a characterization of holomorphic complex variable positive definite functions on a horizontal strip.
The next definition will be useful for that purpose.

\begin{defn}
\label{def_FL_transform}
Let $\mu$ be a non-negative measure. We define the complex variable function
\[ f(z) = \int_{-\infty}^{+\infty} e^{izt} \, d\mu(t) \]
to be the {\em Fourier-Laplace transform\/} of $\mu$,
 and denote it by $f = \mathcal{FL}(\mu)$. 
\end{defn}

\begin{thm}
\label{thm_characterization1}
A function $f$ defined on a horizontal strip $S_{a,b}$ is a holomorphic positive definite function if and only if it is the Fourier-Laplace transform of 
an exponentially finite measure with respect to the interval $I = (a,b)$. Furthermore, the measure $\mu$ is uniquely determined by $f$.
\end{thm}

%
%

\subsection{Examples}
\label{subsec_appl}

In order to show how theorem \ref{thm_characterization1}  applies to positive and negative definite functions, 
we now revisit some of the examples in section \ref{sec_complex_PDF_coPDF}.

\begin{ex}
\label{ex_simplest}
Consider the case $n=1$ in example (1), i.e., the holomorphic function 
\begin{equation}
\label{eq_1+z2}
f(z) = \frac{1}{1+z^2} = \frac{1}{(1+iz)(1-iz)} = \frac{1}{2} \left(\frac{1}{1-iz}+ \frac{1}{1+iz} \right).
\end{equation}

Define $f_c(z) = - \frac{1}{c+iz}$ for $\Im(z) \neq c$. Define also

\[ \sigma^r_c(t) = \left\{
                \begin{array}{ll}
                  0 & \mbox{ if $t<0$} \\
                  e^{ct} & \mbox{ if $t \geq 0$}
                \end{array}
              \right. 
						, \ \ \ \ \
	\sigma^l_c(t) = \left\{
                \begin{array}{ll}
                  e^{ct} & \mbox{ if $t<0$} \\
                  0 & \mbox{ if $t \geq 0$}
                \end{array}
              \right. 
, \]
and consider the non-negative Lebesgue-Stieltjes measures defined by
\[ \mu_c^r(]a,b]) = \int_a^b \sigma^r_c(t) \, dt , \ \ \ \ \mu_c^l(]a,b]) = \int_a^b \sigma^l_c(t) \, dt.\]
It is easily shown that  $\mu_c^r$ (resp. $\mu_c^l$) is exponentially finite with respect to $]c, +\infty[$
(resp. $]-\infty, c[$). The Fourier-Laplace transforms of these measures may be used to represent $f_c(z)$ in the following way:
\begin{equation}
\label{eq_FLT}
 f_c(z) = \frac{1}{c+iz} = \left\{
                \begin{array}{ll}
                  -\int_{-\infty}^{+\infty} e^{izt} \, d\mu_c^r(t) = -\mathcal{FL}(\mu_c^r)  \ \ \mbox{ if $\Im(z) > c$} \\
                   \int_{-\infty}^{+\infty} e^{izt} \, d\mu_c^l(t) = \mathcal{FL}(\mu_c^l)  \ \ \mbox{ if $\Im(z) < c$}
                \end{array}
              \right. 
\end{equation}

Hence $f_c(z)$ is negative definite on $S_{c, +\infty}$ and 
positive definite on $S_{ -\infty, c}$. Observe that $-f_c(z)$ has the exact opposite behaviour on 
these semiplanes.

Now, according to identity \eqref{eq_1+z2}, we may write $f(z) = \frac{1}{2} (f_{1}(z) - f_{-1}(z))$. By choosing the 
appropriate equalities and values of $c$ in \eqref{eq_FLT}, we obtain the following representation for $f$.

\begin{enumerate}
\item[i)] if $ z \in S_{1, +\infty}$, 
\[ f(z)   = \frac{1}{2} \left(-\mathcal{FL}(\mu_{1}^r) + \mathcal{FL}(\mu_{-1}^r) \right) = \frac{1}{2} \mathcal{FL} (-\mu_{1}^r + \mu_{-1}^r). \]
Writing $\mu_{1,+\infty} = \frac{1}{2} (-\mu_{1}^r + \mu_{-1}^r)$ and $\sigma_{1,+\infty} = \frac{1}{2} (-\sigma_{1}^r + \sigma_{-1}^r)$, 
we have that $f(z) = \mathcal{FL}(\mu_{1,+\infty})$, where $$\sigma_{1,+\infty} = \left\{
                \begin{array}{ll}
                  0 & \mbox{ if $t<0$} \\
                  -\sinh(t) & \mbox{ if $t \geq 0$}
                \end{array} \right.  $$
is the density function associated with $\mu_{1,+\infty}$, and we conclude that $f$ is negative definite in $S_{1, +\infty}$.

\item[ii)] if $ z \in S_{-1, 1}$,
\[
 f(z)  = \frac{1}{2} \left(\mathcal{FL}(\mu_{1}^l) + \mathcal{FL}(\mu_{-1}^r) \right) = \frac{1}{2} \mathcal{FL} (\mu_{1}^l + \mu_{-1}^r). \]
Writing $\mu_{-1,1} = \frac{1}{2} (\mu_{1}^l + \mu_{-1}^r)$ and $\sigma_{-1,1} = \frac{1}{2} (\sigma_{1}^l + \sigma_{-1}^r)$, 
we have  $f(z) = \mathcal{FL}(\mu_{-1,1})$, where $$\sigma_{-1,1} =  \left\{
                \begin{array}{ll}
                  e^t &  \mbox{ if $t<0$} \\
                  e^{-t} &   \mbox{ if $t \geq 0$}
                \end{array}  \right. $$
is the density function associated with $\mu_{-1,1}$, and we conclude that $f$ is positive definite in $S_{-1,1}$.

\item[iii)] if $ z \in S_{-\infty, -1}$,
\[ f(z)   = \frac{1}{2} \left(\mathcal{FL}(\mu_{1}^l) - \mathcal{FL}(\mu_{-1}^l) \right) = \frac{1}{2} \mathcal{FL} (\mu_{1}^l - \mu_{-1}^l). \]
Writing $\mu_{-\infty,-1} = \frac{1}{2} (\mu_{1}^l - \mu_{-1}^l)$ and $\sigma_{-\infty, -1} = \frac{1}{2} (\sigma_{1}^l - \sigma_{-1}^l)$, 
we now have  $f(z) = \mathcal{FL}(\mu_{-\infty, -1})$, where 
$$\sigma_{-\infty,-1} =  \left\{
                \begin{array}{ll}
                  \sinh(t) & \mbox{ if $t<0$} \\
                  0 & \mbox{ if $t \geq 0$}
                \end{array} \right.  $$
is the density function associated with $\mu_{-\infty, -1}$, and we conclude that $f$ is negative definite in $S_{-\infty, -1}$.

\end{enumerate}

 Notice that the measures representing $f$ on each of these sets are exponentially 
finite with respect to the intervals  $]1, +\infty[$, $]-1,1[$ and $]-\infty, -1[$, respectively.

\end{ex}

\begin{ex}

Consider the function $f(z) = \frac{1}{\cosh \left( \frac{\pi z}{2} \right)}$ as defined in example 2.
From the identity \eqref{formula_cosh} we may conclude that $f(z) = \mathcal{FL}(\mu_0)$ in $S_{-1,1}$, where $\mu_0$ is
the exponentially finite measure with respect to $]-1,1[$ associated with the density function
$\sigma_0(t) = \frac{1}{\pi \cosh t}$. Observing that, for each $n \in \Z$, $f(z + (2+4n) i) = - f(z)$
and $f(z + 4n i) =  f(z)$, by using the appropriate change of variable we may conclude from 
\eqref{formula_cosh} that:

\begin{enumerate}
\item[i)] $f(z) = \mathcal{FL}(\mu_{2+4n})$ for $z \in S_{1+4n, 3+4n}$ associated with the density function 
$\sigma_{2+4n} (t) = -\frac{e^{(2+4n)t}}{\pi \cosh t}, \ n \in \Z$; 
\item[ii)] $f(z) =\mathcal{FL}(\mu_{4n})$ for $z \in S_{-1+4n, 1+4n}$, where $\mu_{4n}$ is the non-negative exponentially 
finite measure with respect to the interval $]-1+4n, 1+4n[$ 
associated with the density function $\sigma_{4n}(t) = \frac{e^{-4nt}}{\pi \cosh t}, \ n \in \Z$.
\end{enumerate}

\end{ex}

\begin{ex}
Consider the $\Gamma$ function, as defined in example \ref{ex_Gamma}. Notice that, for $\Re(z) > 0$, it may be
written in the form

\[ \Gamma(z) = \int_0^\infty e^{z \, \ln x} \frac{e^{-x}}{x} \, dx. \]

Performing the change of variable $t= - \ln x$ we may write 

\[\Gamma (z) = \int_{+\infty}^{-\infty} e^{-zt} \frac{e^{-e^{-t}}}{e^{-t}} \, (-e^{-t}) \, dt 
							= \int_{-\infty}^{+\infty} e^{-zt} \, e^{-e^{-t}} \, dt. \]

Observe that $\Gamma(-iz) = \int_{-\infty}^{+\infty} e^{izt} \, e^{-e^{-t}} \, dt$  
is, according to theorem \ref{thm_characterization1}, the positive definite 
holomorphic function on $S_{0, +\infty}$ associated with the measure $\mu_0$ such that 
$d\mu_0 (t) = e^{-e^{-t}} \, dt$. Then, in view of proposition \ref{prop_codiff_cosum} , $\Gamma(z)$ is co-positive 
definite on $T_{0,+\infty}$.

For each integer $n \geq 0$ and each vertical strip $T_{-n-1, -n}$ 
we may rewrite the expression in example \ref{ex_Gamma} in the form 

\[ \Gamma(z) = \int_0^\infty e^{z \ln x} \frac{\left( \sum_{m=n+1}^\infty \frac{(-1)^m}{m!} x^m \right)}{x} \, dx \]

Using once again the change of variables $t = -\ln x$ we obtain, on each vertical strip $T_{-n-1, -n}$,

\[ \Gamma(z) = \int_{-\infty}^{+\infty} e^{-zt} \, \left( \sum_{m=n+1}^\infty \frac{(-1)^m}{m!} e^{-mt} \right) \, dt. \] 

Reasoning as in the first case, we conclude that $\Gamma(z)$ is a co-positive or co-negative definite function
on each of the vertical strips $T_{-n-1, -n}$ 
associated, for $n \geq 0$, 
with measures $\mu_{n+1}$ such that $d\mu_{n+1}(t) \! = \!\left( \sum_{m=n+1}^\infty \frac{(-1)^m}{m!} e^{-mt} \right) \! dt.$
Each of these measures is exponentially finite with respect to the interval $]-n-1, -n[$. Their respective signs are determined 
by the sign of the factor $(-1)^{n+1}$ corresponding to the first term of the $n$th order tail of the series  
expansion of $e^x$.
\end{ex}

\begin{ex} 
\label{ex_Bessel_1st_kind}
Consider 
the normalized Bessel function $\mathcal{J}_\alpha$
defined in example \ref{ex_Bessel}. Observe that, for $z \in \C$, its integral representation
may be written in the form 
\begin{align*}
 \mathcal{J}_\alpha(z) & = \int_{-1}^1 \frac{1}{2^{\alpha} \sqrt{\pi} \Gamma(\alpha+1/2)}  (1-t^2)^{\alpha-1/2} \, e^{izt} \, dt \\
						 & = \int_{-\infty}^{+\infty} \sigma_{\mathcal{J}_\alpha}(t) \, e^{izt} \, dt
\end{align*} 
where 
\[ \sigma_{\mathcal{J}_\alpha}(t) =  \left\{
                \begin{array}{ll}
                  \frac{(1-t^2)^{\alpha-1/2}}{2^{\alpha} \sqrt{\pi} \Gamma(\alpha+1/2)} & \mbox{ if $-1 \leq t \leq 1$} \\
                  0 & \mbox{ otherwise.}
                \end{array} \right.  
\]
Notice that $\sigma_{\mathcal{J}_\alpha}$ is the density function of a non-negative measure $\mu_{\mathcal{J}_\alpha}$ Since $\sigma_{\mathcal{J}_\alpha}$
has compact support, $\mu_{\mathcal{J}_\alpha}$ is exponentially finite with respect to the interval $]-\infty, +\infty[$ and the 
corresponding Fourier-Laplace integral
\[  \mathcal{J}_\alpha(z) = \int_{-\infty}^{+\infty} e^{izt} d\mu_{\mathcal{J}_\alpha}(t) \]
represents, by theorem \ref{thm_characterization1}, an entire positive definite function. Since, in the special
case $\alpha = 0$, we have $J_0(z) \equiv \mathcal{J}_0(z)$, we derive the same conclusion for the 
Bessel function $J_0(z)$. 
Notice that, in this case, we have
\[ J_0(z) = \mathcal{FL}(\mu_{J_0})\] 
  where $\mu_{J_0}$ admits the density function
\[ \sigma_{J_0}(t) =  \left\{
                \begin{array}{ll}
                  \frac{1}{\pi \sqrt{1-t^2}} & \mbox{ if $-1 \leq t \leq 1$} \\
                  0 & \mbox{ otherwise.}
                \end{array} \right.  
\]


\end{ex}

\begin{ex}
\label{ex_Zeta_2}
Recall from example \ref{ex_Zeta} that the Riemann $\zeta$ function is co-positive definite on the half-plane $T_{1, +\infty}$,
where it is defined by $\zeta(z) = \sum_{n \geq 1} \frac{1}{n^z}$. Hence, according to theorem \ref{thm_characterization1}
and proposition \ref{prop_codiff_cosum}, the 
$\zeta$ function is associated with a non-negative measure $\mu$ which is exponentially finite with respect to the interval
$]1, +\infty[$. We can identify this unique measure by defining 
$h(t) = \left\{
                \begin{array}{ll}
                  0 \   \mbox{ if $t<0$} \\
                  1 \   \mbox{ if $t \geq 0$}
              		\end{array}
				\right. $
and $H(t) = \sum_{n=1} ^\infty h(t- \ln n)$ since, by writing $\mu(]a,b]) = H(b)-H(a)$ and observing that $h^\prime (t)= \delta (t)$
in the sense of distributions, we have: 

\begin{align*}
\mathcal{FL}(\mu) (iz) & = \int_{-\infty}^{+\infty} e^{-zt} \, d\mu(t) = \int_{-\infty}^{+\infty} e^{-zt} \, H^\prime (t) \, dt  \\
							& = \int_{-\infty}^{+\infty} e^{-zt} \, \sum_{n=1}^\infty \delta(t- \ln n)\, dt = \int_{-\infty}^{+\infty}  \sum_{n=1}^\infty e^{-zt} \, \delta(t- \ln n)\, dt\\
							& = \sum_{n=1}^\infty e^{-z \ln n}  \\
							& = \zeta(z).
\end{align*}

We thus conclude that $H(t)$ is the distribution function for the measure $\mu$ whose corresponding density function
is $\sigma(t) =  \sum_{n=1}^\infty \delta(t- \ln n)$.

\end{ex}

\begin{ex}
\label{ex_ZED}
We 
now consider an interesting analytic extension of the integral representation of the $\zeta$ function
obtained in example \ref{ex_Zeta_2} to the half-plane $\mathcal{R}(z) > 0$. Integrating by parts we have, 
for $\mathcal{R}(z) > 1$, 
\begin{align} 
\label{eq_zeta_computation}
\zeta(z) & = \int_{-\infty}^{+\infty} H^\prime(t) e^{-zt} dt  \\
						& = \left[ H(t) e^{-zt} \right]_{-\infty}^{+\infty} - \int_{-\infty}^{+\infty} H(t) (-z) e^{-zt} dt \\
						& = z \int_0^{+\infty} H(t) e^{-zt} dt.
\end{align}

Observing that $z \int_0^{+\infty} e^{(1-z)t} dt = \frac{z}{z-1}$ for $\mathcal{R}(z) > 1$, we derive 
from \eqref{eq_zeta_computation} that

\[ \zeta(z) = \frac{z}{z-1} - z \int_0^{+\infty} (e^t - H(t))e^{-zt} dt. \]

We denote as usual by $\{.\}$ the fractional part of a real number and 
observe that $e^t - H(t) = e^t - [e^t] = \{e^t\}$ is the fractional part of $e^t$. 
Defining
\[ \sigma_{\mathcal{U}}(t) = \left\{
                \begin{array}{ll}
                  \{e^t\} \   \mbox{ if $t \geq 0$} \\
                  0 \   \mbox{ otherwise}
              		\end{array}
				\right. \]
we have
\begin{equation}
\label{eq_ZEE}
 \zeta(z) = \frac{z}{z-1} - z \mathcal{U}(z)
\end{equation}

\[ \mathcal{U} (z) = \int_{-\infty}^{+\infty} \sigma_{\mathcal{U}}(t) e^{-zt} dt.  \]
Since $\{e^t	\}$ is bounded,  we conclude that $\sigma_{\mathcal{U}}(t)$ is the density function associated with a 
non-negative measure $\mu_{\mathcal{U}}$    which is exponentially finite
with respect to the interval $]0,+\infty[$. Hence $\mathcal{U} (z) \equiv \mathcal{FL}(\mu_{\mathcal{U}})$ is a 
holomorphic co-positive definite function on the half-plane 
$\mathcal{R}(z) >0$ and, by definition of $\zeta(z)$, the meromorphic right hand side of \eqref{eq_extension} coincides, on this set, 
with $\zeta(z)$.

\end{ex}

\begin{ex}
\label{example_ZED}
Recall
 the function $\mathcal{Z} = -\zeta(z) / z$ defined by \eqref{eq_def_ZED}  in example \ref{ex_function_ZED}.
According to \eqref{eq_ZEE}, 
we have
\begin{equation}
\label{eq_extension}
\mathcal{Z}(z) = \frac{1}{1-z} + \mathcal{U} (z) \ \ \ \mbox{ for  $\mathcal{R}(z) > 0$.}
\end{equation}
Recalling the definitions of $\sigma_1^l$ and $\sigma_1^r$ in example \ref{ex_simplest}, and replacing
$z$ with $iz$ in formula \eqref{eq_FLT}, we obtain

\[ \frac{1}{1-z} = \int_{-\infty}^{+\infty} \sigma_1^l (t) e^{-zt} dt \ \ \ \mbox{ for  $\mathcal{R}(z) < 1$}, \]
\[ \frac{1}{1-z} = - \int_{-\infty}^{+\infty} \sigma_1^r (t) e^{-zt} dt \ \ \ \mbox{ for  $\mathcal{R}(z) > 1$}. \]

From \eqref{eq_ZEE} we derive:
\[ \mathcal{Z}(z) = \int_{-\infty}^{+\infty} \left( \sigma_1^l(t) + \sigma_{\mathcal{U}} (t) \right) e^{-zt} dt = 
\int_{-\infty}^{+\infty} \sigma^l_{\mathcal{Z}} (t)  e^{-zt} dt \ \ \mbox{ for  $0< \mathcal{R}(z) <1$} \]
where $\sigma^l_{\mathcal{Z}} (t) = e^t -H(t) = \{e^t\}$, and 

\[ \mathcal{Z}(z) = \int_{-\infty}^{+\infty} \left( -\sigma_1^r(t) + \sigma_{\mathcal{U}} (t) \right) e^{-zt} dt = 
\int_{-\infty}^{+\infty} \sigma^r_{\mathcal{Z}} (t)  e^{-zt} dt \ \ \mbox{ for  $ \mathcal{R}(z) >1$,} \]
where $\sigma^r_{\mathcal{Z}} (t) = -H(t)$.

Hence we conclude that $\mathcal{Z}(z)$ is a  holomorphic co-positive definite function on the vertical strip 
$0< \mathcal{R}(z) <1$ and a holomorphic co-definite negative function on the half-plane $ \mathcal{R}(z) >1$.
Accordingly, the non-negative measure $\mu^l_{\mathcal{Z}}$ (resp. $\mu^r_{\mathcal{Z}}$) associated with the 
density function $\sigma^l_{\mathcal{Z}}$ (resp. $\sigma^r_{\mathcal{Z}}$) is exponentially finite with
respect to the interval $]0,1[$ (resp. $]1, +\infty[)$. 

Notice the interesting fact that the zeros 
of $\zeta(z)$ and $\mathcal{Z}(z)$ coincide. Hence we conclude, 
using the arguments in example \ref{ex_Zeta} that $\mathcal{Z}$ is neither co-positive nor co-negative definite 
in any cosum set contained in the half-plane $ \mathcal{R}(z) < 0$.
\end{ex}

\begin{rem}
\label{rem_IR_Zeta}
We remark that example \ref{example_ZED} proves the following integral representation 
for the $\zeta$ function on the critical strip:
\begin{equation}
\label{eq_IR_Zeta}
\zeta(z) = -z \int_{-\infty}^{+\infty} \{e^t\} e^{-zt} dt, \ \ \ 0 < \Re(z) < 1. 
\end{equation}
A change of variable leads to the integral representation
\begin{equation}
\label{eq_IR_Zeta2}
\zeta(z) = -z \int_0^{+\infty} t^{z-1} \{1/t\}  dt, \ \ \ 0 < \Re(z) < 1,
\end{equation}
\end{rem}
a more usual of several equivalent forms in which it appears in the literature (see e.g. Titchmarsh \cite{Tit}, pg. 15).

We close this section with two results on regularity extension whose proofs may be set in the context of integral representations; 
see theorems 2.11 and 2.13 in \cite{bp_CPDFS}. A different, constructive approach via holomorphic reproducing kernels may be found in 
\cite{bp_PROP} but will not be addressed in this paper.

In the case of real variable positive definite functions, analyticity in a neighborhood of the origin is known to imply analyticity in $\R$
(see Theorem \ref{thm_sasvari}). The next result may be seen as a complex variable counterpart of this property.

\begin{thm}
\label{thm_characterization_2}
Let $f: S_{a,b} \to \C$ be a positive definite function on the horizontal strip $S_{a,b}$. If $f$ is holomorphic 
on some open neighborhood  $\Omega$ of $\Gamma = \{ iy \in \C : a < y < b \}$, 
 then $f$ is holomorphic in $S_{a,b}$ and is the Fourier-Laplace transform of a unique exponentially finite measure $\mu$.
\end{thm}

Consequences of the results in the previous sections
for the properties of meromorphic functions under a minimal condition of positive definiteness are the concern of our next result.

\begin{thm}
\label{thm_minimal_PD}
Let $f$ be meromorphic in $\C$ and suppose $f(iy) = G(y)$ is co-positive definite for $y \in I = (a, b)$, with $a, b \in \overline{\R}$. Then $f$ is a holomorphic positive definite function on the maximal horizontal strip $S_{c,d}$ containing $S_{a,b}$ such that $ic$ and $id$ are poles of $f$ whenever $c$ or $d$ are finite, respectively.
\end{thm}

%
%
%
%

\section{Moments, moment-generating functions and characteristic functions}
\label{sec_moments_MGF_CF}

Suppose $f$ is holomorphic and positive definite on the horizontal strip $S_{a,b}.$ Then, according to Theorem \ref{thm_characterization1}, we have
\begin{equation}
\label{eq_irrep}
f(z) = \int_{-\infty}^{+\infty} e^{izt} \, d\mu(t)
\end{equation}
for a unique 
non-negative exponentially  finite measure with respect to $I = (a,b)$.
 For each $y \in I$, define $M_n^y = \int_{-\infty}^{+\infty} t^n \, d\mu_y(t)$, where
$\mu_y$ 
is the finite measure defined by
\begin{equation}
\label{eq_measure}
\mu_y(E) = \int_E e^{-yt} \, d\mu(t).
\end{equation}
We call $M_n^y$ the $n^{th}$ moment of the measure $\mu_y$. Notice that $f^{(n)}(z) = i^n \, \int_{-\infty}^{+\infty} t^n e^{izt} \, d\mu(t)$ and
$f^{(n)}(iy) = i^n \, \int_{-\infty}^{+\infty} t^n e^{-yt} \, d\mu(t) = i^n M_n^y$. Observe that, by
Theorem \ref{thm_complex_equivalence},  these integrals always exist for $y \in I$.

Defining $h_y(x) = f(z)$ for $z = x+iy$ and $y \in (a,b)$, we may also write
\begin{equation}
\label{eq_h}
h_y(x) = F_y(x) = \int_{-\infty}^{+\infty} e^{ixt} e^{-yt} d\mu(t) = \int_{-\infty}^{+\infty} e^{ixt} \, d\mu_y(t).
\end{equation}
The function $h_y$ is the so-called {\em characteristic function\/} 
for the finite measure $\mu_y$.
Notice that 
\[ h^{(n)}_y(x) = i^n \int_{-\infty}^{+\infty} t^n e^{ixt} \, d\mu_y(t) \]
and consequently
\[ h^{(n)}_y(0) = i^n \int_{-\infty}^{+\infty} t^n  \, d\mu_y(t) = i^n M_n^y.\]
Finally, define $g_y(u) = G(y-u) = f(i(y-u))$ for $y \in (a,b)$ and $y-u \in (a,b)$. Then
\begin{equation}
\label{eq_gyu}
g_y(u) = \int_{-\infty}^{+\infty} e^{(-y+u)t}  \, d\mu(t)  = \int_{-\infty}^{+\infty} e^{ut}  \, d\mu_ y(t).
\end{equation}
		Now 
\[ g^{(n)}_y(u) = \int_{-\infty}^{+\infty} t^n e^{ut}  \, d\mu_y(t) \]
and therefore 
\[  g^{(n)}_y(0) = \int_{-\infty}^{+\infty} t^n   \, d\mu_y(t) = M_n^y. \]
The function $g_y$ is the so-called {\em moment generating function\/} 
for the finite measure $\mu_y$. Hence we may generally write
\[M_n^y = i^n f^{(n)} (iy) = i^n h_y^{(n)}(0) = g_y^{(n)}(0). \]


A special case arises when $y= 0 \in (a,b) = I$. Writing $\mu=\mu_0$ and $h= h_0$
\eqref{eq_h} becomes the usual definition of the 
characteristic function for the finite measure $\mu$, $h(x) = \int_{-\infty}^{+\infty} e^{ixt}  d\mu(t)$. 
From \eqref{eq_gyu}, on the other hand, we obtain $g(u) = g_0(u) = \int_{-\infty}^{+\infty} e^{ut} d\mu(t)$, 
and $g(u)$ identifies with the moment 
generating function for $\mu$. 
In the case where $\mu$ is a probability measure (i.e. a normalized finite measure) we recover the standard context in which these notions have originally been defined.
We may use the previous discussions on complex variable positive definite functions to shed new light on the questions of existence and regularity of the
characteristic and moment generating functions associated with a finite measure $\mu$.

\begin{thm}
\label{thm_moments_cf_measures}
Suppose $\mu$ is a non-negative measure. Then the following statements are equivalent:
\begin{enumerate}
\item $\mu$ is exponentially finite with respect to some open interval containing the origin.
\item $\mu$ admits a characteristic function $h$ which is analytic at the origin.
\item $\mu$ admits a moment generating function $g$ which is defined on some open interval containing the origin.
\end{enumerate}
\end{thm}

\begin{ex}
Recall from examples \ref{ex_Zeta} and \ref{ex_Zeta_2} that the Riemann Zeta function is co-positive
definite on the semiplane $T_{1, +\infty}$ with the integral representation $\zeta(z) = \mathcal{FL}(\mu)(iz)$,
where $d\mu (t) = \sigma(t) \, dt$ and $\sigma(t) = \sum_{n=1}^\infty \delta(t - \ln n)$. Hence 
$f(z) = \zeta(-iz) = \sum_{n=1}^\infty n^{iz}$ is the positive definite  function defined on $S_{1, +\infty}$
associated with $\mu$. According to the above definitions, the density function for the finite measure $\mu_y$ is
\[\sigma_y(t) = \sum_{n=1}^\infty e^{-yt} \delta(t- \ln n) = \sum_{n=1}^\infty n^{-y} \delta(t- \ln n). \]
The characteristic function is $h_y(x) = f(x+iy) = \sum_{n=1}^\infty n^{ix-y}$ and the moment generating function is 
$g_y(u) = f(i(y-u)) = \sum_{n=1}^\infty n^{u-y}.$ For the moments of order $k$ we find the expression
\begin{align}
M^y_k & = \int_{-\infty}^{+\infty} t^k \, d\mu_y(t) = \int_{-\infty}^{+\infty} t^k \sum_{n=1}^\infty n^{-y} \delta(t- \ln n) \, dt \nonumber \\
			& = \sum_{n=1}^\infty n^{-y} \ln^k n. 
\end{align}

\end{ex}


\begin{thebibliography}{99}




%


\bibitem{ask}
R. Askey, 
\newblock Gr\"{u}nbaum's inequality for Bessel functions.
\newblock J. Math. Anal. Appl. 293 (1973), 122--124.

\bibitem{bar}
A. Baricz,
\newblock {\em Generalized Bessel functions of the first kind.}
\newblock Springer-Verlag, LNM {\bf 1994}, Berlin, 2010.



\bibitem{ber}
S. Bernstein, 
\newblock Sur les fonctions absolument monotones.
\newblock Acta Math. 52 (1929), 1--66.


\bibitem{sas00}
T. Bisgaard and Z. Sasvári,
\newblock{\em Characteristic functions and moment problems.}
\newblock Nova Science Publishing, New York, 2000.

\bibitem{bp_DIE}
J. Buescu, A. Paix\~{a}o,
\newblock Positive definite matrices and integral equations on unbounded domains. 
\newblock Differential Integral Equations 19 (2006), no. 2, 189--210. 

\bibitem{bp_JIEA}
J. Buescu, A. Paix\~{a}o, F. Garcia, I. Lourtie,
\newblock Positive-definiteness, integral equations and Fourier transforms. 
\newblock J. Integral Equations Appl. 16 (2004), no. 1, 33--52. 

\bibitem{bp_LAA}
J. Buescu, A. Paix\~{a}o,
\newblock A linear algebraic approach to holomorphic reproducing kernels in ${\mathbb C}^n$.
\newblock  Linear Algebra Appl.  412  (2006),  no. 2-3, 270--290.

\bibitem{bp_pdmdrki}
J. Buescu, A. Paix\~{a}o,
\newblock Positive definite matrices and differentiable reproducing kernel inequalities.
\newblock J. Math. Anal. Appl. {\bf 320} (2006), 279--292.

\bibitem{bp_diff_PDF}
J. Buescu, A. Paix\~{a}o,
\newblock On differentiability and analyticity of positive definite functions. 
\newblock J. Math. Anal. Appl. 375 (2011), no. 1, 336--341. 

\bibitem{bp_RC_PDF}
J. Buescu, A. Paix\~{a}o,
\newblock Real and complex variable positive definite functions. 
\newblock S\~{a}o Paulo Journal of Mathematical Sciences {\bf 6}, 2 (2012), 155--169.

\bibitem{bp_CVPDF}
J. Buescu, A. Paix\~{a}o,
\newblock Complex variable positive definite functions. 
\newblock  Complex Anal. Oper. Theory {\bf 8} (2014), no. 4, 937--954. 

\bibitem{bp_CPDFS}
J. Buescu, A. Paix\~{a}o, A. Symeonides,
\newblock Complex positive definite functions on strips.
\newblock Complex Anal. Oper. Theory 11 (2017), no. 3, 627--649.

\bibitem{bp_PROP}
J. Buescu, A. Paix\~{a}o, C. Oliveira,
\newblock Propagation of regularity and positive definiteness: a constructive approach.
\newblock To appear.



 \bibitem{dev2}
A. Devinatz,
 \newblock Integral representations of positive definite functions.
\newblock   Trans. Amer. Math. Soc. {\bf 74} (1953), 56--77. 

\bibitem{dev3}
A. Devinatz,
\newblock Integral representations of positive definite functions II. 
\newblock  Trans. Amer. Math. Soc. {\bf 77}, (1954), 455--480.


\bibitem{don}
W. Donoghue,
\newblock {\em Distributions and Fourier transforms}.
\newblock Academic Press, New York, 1969.



%

\bibitem{krein}
M. Krein,
\newblock Sur le problème du prolongement des fonctions hermitiennes positives et continues. 
\newblock C. R. (Doklady) Acad. Sci. URSS (N.S.) 26 (1940). 17--22. 

\bibitem{kos}
H. Kosaki (2011),
 \newblock {\em  Positive definiteness of functions with applications to Operator Norm Inequalities.}
 \newblock Mem. Amer. Math. Soc., 212.



\bibitem{luk2}
E. Lukacs,
\newblock{\em Characteristic functions.}
\newblock{Griffin and Co., 2nd ed., 1970}.

\bibitem{mat}
M. Mathias,
\newblock Über positive Fourier-Integrale.
\newblock Math. Zeitsch. {\bf 16} (1923), 1,  103--125.

\bibitem{neu}
E. Neuman, 
\newblock{Inequalities involving Bessel functions of the first kind.}
\newblock J. Inequal. Pure Appl. Mat 5(4), article 94, 4pgs (electronic). 

\bibitem{ste}
J. Stewart, 
\newblock Positive definite functions and generalizations, an historical survey.
\newblock Rocky Mountain Journal of Mathematics {\bf 6} (1976), no. 3, 409--434.





\bibitem{sas94}
Z. Sasvári,
\newblock{\em Positive definite and definitizable functions.}
\newblock  Mathematical Topics, 2. Akademie Verlag, Berlin, 1994. 


\bibitem{sas13}
Z. Sasvári,
\newblock{\em Multivariate characteristic and correlation functions.}
\newblock De Gruyter Studies in Mathematics, 50. Walter de Gruyter \& Co., Berlin, 2013. 

\bibitem{sel}
V. Selinger, 
\newblock{Geometric properties of normalized Bessel functions.}
\newblock Pure Math. Appl. 6 (1995), 273--277. 



\bibitem{Tit}
E. C. Titchmarsh,
\newblock {\em The Theory of the Riemann Zeta-Function.}
\newblock The Clarendon Press, Oxford University Press (2nd ed), 1986.
%

\bibitem{wat}
G. Watson, 
\newblock {\em A treatise on the theory of Bessel functions.}
\newblock Cambridge U. P., Cambridge, 1944.

\bibitem{whi-wat}
E. Whittaker, G. Watson, 
\newblock {\em A course in modern analysis.}
\newblock Cambridge U. P., Cambridge (4th ed), 1996.

\bibitem{wid} 
D. Widder,
\newblock Necessary and sufficient conditions for the representation of a function by a doubly infinite Laplace integral.
\newblock Bull. Amer. Math. Soc. {\bf 40} (1934), no. 4, 321--326.




\end{thebibliography}
\end{document}